\documentclass[11pt]{amsart}
\usepackage{amssymb}
\theoremstyle{plain}

\newcommand{\N}{\mathbb{N}}

\newcommand{\R}{\mathbb{R}}

\newcommand{\Id}{\text{\rm Id}}

\newcommand{\cF}{\mathcal{F}}

\newcommand{\cH}{\mathcal{H}}

\newtheorem{thm}{Theorem}[section]

\newtheorem{cor}[thm]{Corollary}
\newtheorem{prop}[thm]{Proposition}
\theoremstyle{definition}

\newtheorem*{rem}{Remark}
\newtheorem*{rems}{Remarks}

\newcommand{\ra}{\rangle}
\newcommand{\la}{\langle}

\newcommand{\supp}{{\rm supp}}

\newcommand{\kin}{\!\in\!}

\def\hangbox to #1 #2{\vskip1pt\hangindent #1\noindent \hbox to #1{#2}$\!\!$}

\allowdisplaybreaks

\title[nonuniform sampling and recovery via Gaussian RBFs]{Nonuniform
sampling and recovery of\\ multidimensional bandlimited functions 
by\\ Gaussian radial-basis functions}

\author{B.~A.~Bailey}\author{Th. Schlumprecht}\author{N.~Sivakumar}
\address{Center for Approximation Theory, Department of Mathematics, Texas A\&M University\\  
College Station, TX 77843, USA}
\email{abailey@math.tamu.edu, schlump@math.tamu.edu, sivan@math.tamu.edu}

\thanks{\textit{2000 Mathematics Subject Classification}: Primary 41A05, Secondary 42C30} 
\thanks{The research of
the  second author was supported by the National Science Foundation}
\keywords{Scattered Data, Gaussian interpolation, Multidimensional Bandlimited functions} 

\begin{document}
\maketitle
\begin{abstract}
Let $S\subset\R^d$ be a bounded subset with positive Lebesgue measure.
The Paley-Wiener space associated to $S$, $PW_S$, is defined
to be the set of all
square-integrable functions on $\R^d$ whose Fourier transforms vanish
outside $S$.
A sequence $(x_j:j\kin\N)$ in $\R^d$ is said to be a Riesz-basis sequence for 
$L_2(S)$ (equivalently, a complete interpolating sequence for $PW_S$) 
if the sequence $(e^{-i\la x_j,\cdot\ra}:j\kin\N)$
of exponential functions 
 forms a Riesz basis for $L_2(S)$.
Let $(x_j:j\kin\N)$ be a Riesz-basis sequence for 
$L_2(S)$.  Given
$\lambda>0$  and $f\in PW_S$,  
there is a unique
sequence $(a_j)$ in $\ell_2$ such that the function
$$
I_\lambda(f)(x):=\sum_{j\in\N}a_je^{-\lambda \|x-x_j\|_2^2}, \qquad x\kin\R^d,
$$
is continuous and square integrable on $\R^d$,
and satisfies the  condition
$I_\lambda(f)(x_n)=f(x_n)$ for every $n\kin\N$.
This paper studies the convergence
of the interpolant $I_\lambda(f)$ as 
$\lambda$ tends to zero, {\it i.e.,\/} as the variance
of the underlying Gaussian tends to infinity.  The
following result is obtained:
Let $\delta\in(\sqrt{2/3},1]$ and $0<\beta<\sqrt{3\delta^2 -2}$.  Suppose that
$\delta B_2\subset Z\subset B_2$,
and let  $(x_j:j\in\N)$ be a Riesz basis sequence for $L_2(Z)$.  
If $f\in PW_{\beta B_2}$, then  
$f=\lim_{\lambda\to 0^+} I_\lambda(f)$ in $L_2(\R^d)$ and uniformly on $\R^d$.
If $\delta=1$, then one may take $\beta$ to be $1$ as well, 
and this reduces to a known theorem in the univariate case.
However, if $d\ge2$, 
it is not known
whether $L_2(B_2)$ admits a Riesz-basis sequence.  On the other hand,
in the case when $\delta<1$, there do exist bodies $Z$ satisfying the 
hypotheses of the theorem (in any space dimension). 

\end{abstract}

\section{Introduction}\label{S:0}

The theoretical study of interpolation has long played a 
prominent role in the development of the theory
of approximation, whilst the computational aspects of 
the subject have found a natural outlet in numerical analysis.
Given the inherent naturalness of the subject, and the 
variety of interesting questions it generates, it is not
surprising that interpolation theory continues to 
be an active area of research.

Among the diverse issues associated to interpolation, the theory
of `cardinal interpolation' has been a well-studied theme
over many years.  The term refers to 
the interpolation of data on a regular, usually infinite, grid.
When the interpolation at such a grid is done via spline functions,
in particular, one encounters `cardinal spline interpolation'.  
Rooted in 
Isaac Schoenberg's seminal work, this subject 
was further developed by a number of his
successors.

The last century has witnessed enormous
advances in the state of the art of electrical engineering
and telecommunication theory.  Included among the myriad branches
of this discipline is the theory of `sampling', which is an
integral part of signal analysis.  By their very nature,
the theories of sampling and interpolation are intertwined;
indeed, at a basic level both are in essence one and the same.
Philosophical considerations apart, there are also solid
and captivating mathematical connections between the 
subjects.  Connections of this nature, in the context
of cardinal spline interpolation, were brought out
by Schoenberg himself; see, for instance 
\cite{Sc}, especially his 
remarks on page 228 there.  In this paper Schoenberg showed,
among other things, how bandlimited functions
can be recovered as limits of their cardinal-spline
interpolants, as the degree of the underlying
spline tends to infinity.  Substantial extensions
of this theme have since been carried out by many,
both in the univariate and multivariate settings.  More recently,
it was also discovered that the theory of cardinal spline
interpolation, and its connections to sampling, have 
a strong resonance in the theory of radial-basis
functions, especially involving Gaussians.  One 
consequence of these developments is the fact
that bandlimited functions can also be recovered
via their Gaussian cardinal interpolants, in a suitable
limiting sense.

Moving away from the realm of gridded data, it is natural
to look for comparable connections between interpolation
at irregularly spaced points -- often
referred to as scattered-data interpolation -- via
splines and radial-basis functions, and the now classical 
theory
of nonuniform sampling.  In the case of splines,
the analogue of Schoenberg's theorem to which
we alluded above was obtained in \cite{LM}, and
a counterpart of this
result for Gaussians was given later
in \cite{SS}.

The focus in \cite{SS} is on the one-dimensional
case, although it  does include a relatively
straightforward extension
to higher dimensions in terms of tensor products.  
The purpose of 
the present article is to present a general
multidimensional result. Here we overcome our
earlier obstacles by extending the 
techniques of \cite{LM} and 
\cite{SS} in such a way that
the radial symmetry of the multidimensional
Gaussian can be exploited fully.  As in \cite{LM} and \cite{SS}, our setting for the
interpolation problem also  involves Riesz-basis sequences. 
However, we are forced to settle for less in the
multivariate situation, because the existence 
of suitable Riesz-basis sequences in higher dimensions is 
a subtle matter depending on the geometry of the underlying domain.  
A more detailed discussion of this issue is the subject of 
the concluding remarks in Section 3.
  
The rest of the paper is organized as follows.
In the next section we lay out some background
material of a general nature.  More specific
preliminaries are given in the subsequent section, 
which concludes with
the statement of the main result.  The proof
of the latter
is detailed in the final section.

\section{Preliminaries}\label{S:1}

Let $d\in\N$.  For $1\le p< \infty$ and  a measurable set  
$S \subset\R^d$ with positive Lebesgue measure $m(S)$,  we denote by
 $L_p(S)$   the space of complex-valued functions which are $p$-integrable  on $S$ (with respect
 to the Lebesgue measure). For $f\in L_p(S)$, we denote its standard $L_p$-norm
 by 
 $\|f\|_{L_p(S)}$, or by $\|f\|_p$, when the context is clear.
We also denote by $\|\cdot\|_p$, 
the $\ell_p$-norm, $1\le p\le \infty$, on the space of (finite and infinite) sequences.
The space of continuous functions on $\R^d$ is denoted by $C(\R^d)$, and
$$
C_0(\R^d):=\big\{f\in C(\R^d): \lim_{\|x\|_2\to\infty} f(x)=0\big\}.
$$

If $f\in L_1(\R^d)$, then the Fourier transform of $f$,  $\hat f$, is defined as follows:
\begin{equation}\label{E:1.1}
 \hat f(x):=\int_{\R^d} f(u) e^{-i\la x,u\ra}du, \quad \text{ $x\in\R^d$.}
\end{equation}
The Fourier transform  of   $g\in L_2(\R^d)$ is denoted by $\cF[g]$.  
The assignment $g\mapsto\cF[g]$ is the unique extension
 of the map
 $$
\widehat{(\cdot)}: L_1(\R^d)\cap L_2(\R^d)\to  L_2(\R^d)
$$
to a bounded linear  operator on $L_2(\R^d)$, satisfying
the following Plancherel-Parseval relation:
\begin{equation}\label{E:1.2}
\|\cF[g]\|_2^2= (2\pi)^{d}\|g\|_2^2\text{  \ for all $g\in L_2(\R^d)$.}
 \end{equation}
If $g\in L_2(\R^d)\cap C(\R^d)$ and $\cF[g]\in L_1(\R^d)$, then the following inverse formula
 holds:
\begin{equation}\label{E:1.3}
 g(x)=\frac1{(2\pi)^d} \int_{\R^d} \cF[g](u) e^{i\la u,x\ra} du, \quad\text{ for all $x\in\R^d$.}
\end{equation}

For $\lambda>0$ we define the {\em Gaussian function } $g_\lambda:\R^d\to\R$ by
$$
g_\lambda(x)=e^{-\lambda\|x\|_2^2}, \quad \text{for all $x\in\R^d$},
$$
and recall that 
\begin{equation}\label{E:1.4}
 \hat g_\lambda(u)=\Big(\frac{\pi}{\lambda}\Big)^{d/2} e^{-{\|u\|_2^2/(4\lambda)}}, \text{ for all $u\in\R^d$}.
\end{equation}

The functions we wish to interpolate are the so called {\em bandlimited } or
 {\em Paley-Wiener functions} on $\R^d$. Specifically, for a bounded and measurable
   $S\subset \R^d$ with $m(S)>0$, we define
$$
PW_S=\big\{g\in L_2(\R^d): \cF[g]=0 \text{ almost everywhere outside $S$}
\big\}.
$$

If $S$ is as above and $g\kin PW_S$,
then
the Fourier inversion formula implies the relations
\begin{equation}\label{E:1.5}
g(u)=\frac{1}{(2\pi)^d}\int_{\R^d}
{\cF}[g](x)e^{i\la x, u\ra}\,dx=
\frac{1}{(2\pi)^d}\int_A
{\cF}[g](x)e^{i\la x, u\ra}\,dx,
\end{equation}
for almost all $u\kin\R^d$. As $L_2(S)\subset L_1(S)$ and
${\cF}[g]=0$ almost everywhere outside $S$, it follows
that ${\cF}[g]\kin L_1(\R^d)$, so the 
Riemann-Lebesgue Lemma asserts that the last
expression in \eqref{E:1.5}, as a 
function of $u$, belongs to $C_0(\R^d)$.
So we may assume that \eqref{E:1.5} holds for all $u\kin\R^d$.
Moreover,  
the Bunyakovskii--Cauchy--Schwarz (BCS)  
Inequality  and \eqref{E:1.5} 
combine to show that
\begin{equation}\label{E:0.5a}
|g(u)|\le \frac{m^{1/2}(S)}{(2\pi)^{d}}\|\cF[g]\|_{L_2(S)}=
\frac{m^{1/2}(S)}{(2\pi)^{d/2}}\|g\| _{L_2(\R^d)}, 
\quad u\kin\R^d.
\end{equation}

\section{Gaussian interpolants associated to Riesz-basis sequences}\label{S:2}

We begin by assembling some basic facts about bases in Hilbert spaces (cf. \cite{Yo}).
Let $(\cH,\la\cdot,\cdot\ra_H)$ be a separable (complex) infinite dimensional  Hilbert space. We say that $(h_j:j\kin \N)$ is a {\em Riesz basis for $\cH$} if 
every element $h$ in $\cH$
 admits a unique representation of the form
 \begin{equation}\label{E:2.1}
h=\sum_{j\in \N} a_j h_j,   \quad\text{ with }\qquad \sum_{j\in \N} |a_j|^2<\infty.
 \end{equation}
 One can then show that there exists a unique bounded sequence $(h^*_j:j\kin\N)\subset \cH$ so that 
$a_j=\la h,h^*_j \ra_H$, for  all $j\kin\N$.  We call the $h^*_j$'s the {\em coordinate functionals for $(h_j)$}.
 The sequence $(h^*_j)$ is also a Riesz basis for $\cH$,
and its coordinate functionals are  the $h_j$'s.
 Moreover there exists
 a positive  constant $R_b$ so that 
 \begin{equation}\label{E:2.2}
\frac1{R_b}\Big( \sum_{j\in\N}  |c_j|\Big)^{1/2}\le \Big\|\sum_{j\in\N} c_j h_j \Big\|\le R_b \Big( \sum_{j\in\N}  |c_j|\Big)^{1/2}
\end{equation}
 for every  square-summable  sequence $(c_j:j\kin\N)$.  

Let $S$ be a bounded subset of $\R^d$ with positive Lebesgue measure. We say that a sequence $(x_j:j\kin\N)\subset \R^d$
is a {\em Riesz-basis sequence for $L_2(S)$}  if the sequence $(e^{-i\la x_j,\cdot\ra}:j\kin\N)$  
is a Riesz basis for $L_2(S)$.
The following   pair  of observations will be useful.

\begin{prop}\label{P:2.0} Let $S$ be as above and let $(x_k)$ be a Riesz-basis sequence for $L_2(S)$.
\begin{enumerate}
\item{(a)} There is a $q>0$ so that $\|x_k-x_\ell\|\ge q$ for $k\not=\ell$.
\item{(b)} There is a constant $C>0$ so that
$$\Big(\sum_{k\in\N} |f(x_k)|^2\Big)^{1/2}\le C \|f\|_{L_2(\R^d)},
   \text{ for every $f\in PW_{S}$.}$$
\end{enumerate}
\end{prop}
\begin{proof} If (a) were not true, there would be  two subsequences $(k_j)$ and $(\ell_j)$ 
  such that $\lim_{j\to\infty} \|x_{k_j}-x_{\ell_j}\|= 0$,  
whence the Dominated Convergence Theorem   implies that 
 $\lim_{j\to\infty} \|e^{-i\la x_{k_j},\cdot\ra}-e^{-i\la x_{\ell_j},\cdot\ra}\|_{L_2(S)}= 0$. Let 
  $(e^*_j:j\kin\N)$  be the coordinate functionals for $(e^{-i\la x_j,\cdot\ra}:j\kin\N)$. Since $\la  e^{-i\la x_n,\cdot \ra},e^*_m\ra=\delta_{mn} $, for $m,n\in\N$,
 we  have (with $\la\cdot,\cdot\ra_{_S}=\la\cdot,\cdot\ra_{L_2(S)}$)
  $ \la e^{-i\la x_{k_j},\cdot\ra}-e^{-i\la x_{\ell_j},\cdot\ra}, e^*_{\ell_j}\ra_S=1$,  for $j\in\N$.
But this contradicts the boundedness of  the sequence
   $(\|e^*_j\|:j\kin\N)$.

\noindent (b) Let $f\in PW_S$. Considering $\cF[f]$ as a function in $L_2(S)$
 and  recalling that $(e^*_j)$ is also  a Riesz basis for $L_2(S)$,  we may write
$$\cF[f]= \sum_{j\in\N} \la \cF[f], e^{-i\la x_j,\cdot \ra}\ra_{_S} e^*_j= (2\pi)^d\sum_{j\in\N}  f(x_j) e^*_j,$$
where the second equality stems from \eqref{E:1.5}.
The asserted result follows from \eqref{E:2.1} and \eqref{E:2.2} applied to the Riesz basis $(e^*_j)$.
  \end{proof}

We now begin our discussion of Gaussian interpolants associated to to Riesz-basis sequences. As the first step we state
the following result which can be deduced from \cite[Lemma 2.1]{NSW} and the M.~Riesz Convexity
 Theorem.

\begin{prop}\label{P:2.1a} Let $q>0$ and  suppose that $(x_j)$ is a  $q$-separated sequence in
 $\R^d$, {\it i.e.} 
 $$\inf_{k\not=\ell} \|x_k-x_\ell\|\ge q.$$
  If $\lambda>0$, and if $(a_j)$ is a bounded sequence of complex numbers, then the function
 $\R^d\ni x\mapsto\sum a_j g_\lambda(x-x_j)$ is continuous and bounded.
 Moreover, the infinite matrix $\big(g_\lambda(x_j-x_k)\big)_{j,k\in\N}$ acts as a bounded operator
 on $\ell_p$, for all $1\le p\le \infty$.
\end{prop}

This next theorem is an important discovery in the 
quantitative theory
of radial-basis functions.
\begin{thm}\label{T:NW} { \cite[Theorem 2.3]{NW}}
Let $\lambda$ and $q$ be fixed positive numbers.  There exists a number $\theta$, depending only
on $d$, $\lambda$, and $q$, such that
the following holds: if 
$(x_j)$
is any  $q$-separated sequence in $\R^d$, 
then
$
\sum_{j,k}\xi_j{\overline \xi}_k g_\lambda(\|x_j-x_k\|_2)\ge\theta
\sum_j|\xi_j|^2,
$
for every sequence of complex numbers
$(\xi_j)$. 
\end{thm}

\begin{cor}\label{C:2.2}
Suppose that $\lambda$ is a fixed
positive number. Let
$(x_j:j\kin\N)\subset \R^d$  be $q$-separated for some $q>0$.
 Then the matrix 
$(g_\lambda(x_k-x_j))_{k,j\in\N}$  is boundedly invertible on $\ell_2$.
In particular,
given a 
square-summable sequence $(d_k:k\kin\N)$, there
exists a unique square-summable sequence
$(a_j^{(\lambda)}:j\kin\N)$ such that
$$
\sum_{j\in\N}a_j^{(\lambda)}g_\lambda(x_k-x_j)=d_k, \qquad k\kin\N.
$$
\end{cor}

The interpolation
operators, whose study will occupy the rest
of the paper, are introduced in 
the following theorem, which is a 
necessary prelude to the main result; its proof will
be given in the next section.
 
\begin{prop}\label{P:2.3} Let $d\kin\N$ and
 let $ Z  \subset \mathbb{R}^d$ be convex, symmetric about the origin and bounded with $m(Z)>0$.  Let $\lambda$ be  a fixed positive number, 
and let $(x_j:j\kin\N)\subset \R^d$ be a Riesz -basis sequence for $L_2(Z)$.

For any $f\in PW_{Z}$, there exists a unique square-summable 
sequence $(a_j^{(\lambda)}:j\kin\N)$
such that 
\begin{equation}\label{E:2.3.1} 
\sum_{j\in\N}a_j^{(\lambda)}  g_\lambda(x_k-x_j)=f(x_k), \quad k\kin\N.
\end{equation}
The {\em Gaussian Interpolation Operator }
$I_\lambda: PW_{ Z }\to L_2(\R^d)$,
defined by
\begin{equation}\label{E:2.3.1a}
I_\lambda(f)(\cdot)=\sum_{j\in\N} a_j^{(\lambda)}g_\lambda((\cdot)-x_j),
\end{equation}
where $(a_j^{(\lambda)}:j\kin\N)$ satisfies \eqref{E:2.3.1},
is a well-defined, bounded linear operator from $PW_{ Z }$ to  $L_2(\R^d)$.
Moreover, $I_\lambda(f)\in C_0(\R^d)$.
\end{prop}

We now state the main result of the paper.

\begin{thm}\label{T:2.4} Let $\delta\in(\sqrt{2/3},1]$ and $0<\beta<\sqrt{3\delta^2 -2}$. Assume that  $Z\subset \R^d$  is convex and  symmetric about the origin, such that 
$\delta B_2\subset Z\subset B_2 $, and 
and let  $(x_j:j\in\N)$ be a Riesz basis sequence for $L_2(Z)$.  Let $I_\lambda$
 be the associated Gaussian Interpolation Operator.  
  Then for every $f\in PW_{\beta B_2}$ we have  $ f=\lim_{\lambda\to 0^+} I_\lambda(f)$ in $L_2(\R^d)$ and uniformly on $\R^d$.
\end{thm}

\begin{rems} The statement of Theorem \ref{T:2.4} includes the case $\delta=1$, {\it i.e.,\/}  $Z=B_2$.  Moreover, in this case, the proof allows one to take 
$\beta=1$.  However, unless $d=1$ (in which case one obtains Theorems 4.3 and 4.4 in \cite{SS}), 
this result may well be vacuous, for 
it is not known if $L_2(B_2)$ admits a Riesz-basis sequence for $d\ge2$.  
In fact, Kristian Seip and Joaquim Ortega-Cerd\`a have informed us
that the prevailing belief is that, when $d\ge2$, there is {\it no\/} Riesz basis 
for $L_2(B_2)$ consisting of exponentials; the latter
has also demonstrated this to be the case for certain allied spaces \cite{OC}.
This interesting problem is 
closely related to interpolatory properties of the associated Paley-Wiener space.  It is 
also connected to 
Fuglede's  work \cite{Fu} and his conjecture, and to the recent studies reported in \cite{IKT1}, \cite{IKT2}, and \cite{Ta}.

On the other hand, there do exist bodies $Z$ satisfying the hypotheses of Theorem \ref{T:2.4}
(in any space dimension); specifically these are {\sl zonotopes\/}.  Firstly,
given any $\delta\in(0,1)$,  
there exist zonotopes $Z$ such that $\delta B_2\subset Z\subset B_2$; this fact, well known to convex geometers,
may be deduced, for instance, as a consequence
of \cite[Theorem 4.1.10]{Ga}.  As to the existence of Riesz-basis sequences
for $L_2(Z)$, this is proved in \cite{LR}  for $d=2$;  the higher dimensional version of this theorem -- to which \cite{LR} alludes --
was communicated to us by Yuri Lyubarskii (private
correspondence).

It is also known that Riesz-basis sequences exist for $L_2(T^d)$, where $T^d$ is a symmetric cube centred at the origin.  Furthermore, in this case,
one can provide sufficient conditions under which a set of distinct points in $\R^d$ forms a Riesz-basis sequence for 
$L_2(T^d)$; see, for example, \cite{SZ} and \cite{Ba}.  These conditions lead to multivariate generalizations of Kadec's
famous ``1/4-theorem" \cite{Ka}.  However, cubes do not quite serve our purpose here; our argument makes essential use
of the additional flexibility offered by zonotopes.

\end{rems}

\section{Proof of the main result}\label{S:3}

For $m\in\N$ we define a linear bounded operator   $A_m$ on $L_2( Z )$ as follows:
Let $(e^*_k) \subset L_2( Z )$ be the 
 coordinate functionals  for $(e^{-i\la x_k,\cdot\ra}:k\kin\N)$, {\it i.e.,}  for every $h\in L_2( Z )$,
\begin{equation}\label{E:3.2}
 h=\sum_{k\in\N} \la h,e^*_k\ra_{_{ Z }} e^{-i\la\cdot, x_k\ra} =
\sum_{k\in\N}    \int_{ Z } h(\xi) \overline{e^*_k(\xi)} \,d\xi \, e^{-i\la\cdot, x_k\ra}.
\end{equation}

Note that for $a=(a_1,a_2\ldots, a_d)\in \R^d$ we have 
\begin{align}\label{E:3.5}
 \Big\|\sum_{k\in\N} \la h,e^*_k\ra_{_{ Z }} e^{-i\la \cdot, x_k\ra}\Big\|_{L_2(a+ Z )}
 &= \Big\|\sum_{k\in\N} \la h,e^*_k\ra_{_{ Z }} e^{-i\la \cdot+a, x_k\ra}\Big\|_{L_2( Z )} \\
 &= \Big\|\sum_{k\in\N} e^{-i\la a, x_k\ra} \la h,e^*_k\ra_{_{ Z }} e^{-i\la \cdot, x_k\ra}\Big\|_{L_2( Z )}
 \le R_b^2 \|h\|_{L_2( Z )},\notag
\end{align}
where $R_b$ is the constant satisfying  \eqref{E:2.2}.
Thus the following extension $E(h)$ 
 of $h$  is locally square integrable, hence   defined almost everyhwere on $\R^d$.
 \begin{equation}\label{E:3.3}
E(h)(x)=\sum_{k\in\N} \la h,e^*_k\ra_{_{ Z }} e^{-i\la x, x_k\ra},\,\, x\in\R^d . 
\end{equation}
Let $m \in \N$, and define $A_{m}:L_2( Z )\to L_2( Z )$ by
 \begin{equation}\label{E:3.6}
 A_m(h)(\xi)=E(h)(2^m(\xi))\chi_{ Z \setminus (1/2) Z }(\xi)
 \end{equation}
 For $h\in L_2( Z )$ it follows from \eqref{E:3.5}  that 
 \begin{align}\label{E:3.7} 
 \|A_m(h)\|_{L_2( Z )}^2&=\int_{ Z \setminus (1/2) Z } |E(h)(2^m u)|^2 du\\
 &=
 2^{-dm}\int_{2^m Z \setminus 2^{m-1} Z } |E(h)(v)|^2 dv\le 2^{-dm} C^m 
  R_b^4\|h\|_{L_2( Z )}^2,\notag
 \end{align}
 where $C$  is the number of translates of $ Z $ which are needed to cover $2Z$.
 The constant $C$  can be bounded by a number which only depends on $d$, and an induction argument shows that at most $C^m$ translates of $ Z $ are  needed to cover $2^m Z $.
\begin{proof}[Proof of Proposition \ref{P:2.3}]
Let $\lambda>0$. By Proposition \ref{P:2.0} and Corollary \ref{C:2.2},  there is a positive constant $\kappa$ so that,
  for each $f\in PW_{ Z }$, there is a sequence $(a_j^{(\lambda)})\in\ell_2$ 
satisfying \eqref{E:2.3.1} and the estimate
\begin{equation}\label{E:2.3.3}
\|(a_j^{(\lambda)})\|_2\le \kappa\|f\|_2.
\end{equation}
Proposition \ref{P:2.1a} ensures that the function $I_\lambda(f)$, as defined in \eqref{E:2.3.1}, is
 continuous and bounded whenever $f\in PW_{ Z }$.

Next we show that $I_\lambda$ is a bounded operator on $L_2(\R^d)$.
Let $f\in PW_{ Z }$ and let $(a_j^{(\lambda)})\in \ell_2$ be the sequence  given above.
By \eqref{E:2.2}, the function $Q:=\sum_{k\in\N} a_k^{(\lambda)} e^{-i\la \cdot, x_k\ra}$ 
is square integrable on $ Z $,
so \eqref{E:3.5} ensures that 
$\|Q\|_{L_2(a+ Z )}\le R_b^2\|Q\|_{L_2( Z )}$ whenever $a\in \R^d$. 
In particular, $Q$ is locally square integrable, hence locally integrable,
on $\R^d$.  Combining these facts with the exponential
decay of $\hat{g_\lambda}$, we find, via a standard periodization argument, that
the function
$$ 
w: \R^d\ni x\mapsto \Big(\frac{\pi}{\lambda}\Big)^{d/2} e^{-\|x\|^2/(4\lambda)}\sum_{k\in\N} a_k^{(\lambda)} e^{-i\la x, x_k\ra},
$$
belongs to $L_2(\R^d)\cap L_1(\R^d)$.
Moreover, using \eqref{E:2.2}, \eqref{E:3.5}, and  \eqref{E:2.3.3}, 
 we arrive at the estimate 
\begin{equation}\label{E:3west}
\|w\|_{L_2(\R^d)}\le C'\|f\|_{L_2(\R^d)},
\end{equation} 
where $C'$
depends only on $\lambda$ and $R_b$.
As $w$ is in $L_1(\R^d)\cap L_2(\R^d)$ and $I_\lambda(f)$ is continuous, it follows from general principles that
$w$ is the Fourier transform of $I_\lambda(f)$. Thus $I_\lambda(f)\in C_0(\R^d)\cap L_2(\R^d)$ and  
$\|I_\lambda(f)\|_{L_2(\R^d)}\le C' (2\pi)^{-d/2}\|f\|_{L_2(\R^d)}$, by \eqref{E:3west} and \eqref{E:1.2}.
\end{proof}

\begin{proof}[Proof of Theorem \ref{T:2.4}]

Now fix $f\in PW_{ Z }$ and write $I_\lambda(f)$ as
$$I_\lambda(f)(\cdot)=\sum_{j\in\N} a_j^{(\lambda)}g_\lambda((\cdot)-x_j).$$
Recall from the preceding paragraph that
the Fourier transform of $I_\lambda(f)$ is given by
\begin{equation}\label{E:3.1}
 \cF\big[I_\lambda(f)\big](u)=
\Big(\frac{\pi}{\lambda}\Big)^{d/2} e^{-\|u\|^2_2/(4\lambda)} \sum_{j\in\N} a^{(\lambda)}_j e^{-i \la x_j,u\ra}, \text{ $u\in\R^d$}.
\end{equation}

The proof of Theorem \ref{T:2.4} proceeds in three steps.

\noindent
{\bf Step 1.} We claim that there is a constant $D_1<\infty$ and $\lambda_0>0$, only depending on $(x_j)$, so that 
$$\Arrowvert \mathcal{F}[I_\lambda(f)]\Arrowvert_2 \le D_1 e^{(1-\delta^2)/(4\lambda)} \Arrowvert \mathcal{F}(f) \Arrowvert_2,  \quad \lambda \in (0,\lambda_0].$$
We start by defining 
$$
H_\lambda(u)=\Big(\frac{\pi}{\lambda}\Big)^{d/2} \sum_{j\in\N} a^{(\lambda)}_j e^{-i \la x_j,u\ra}= e^{\|u\|_2^2/(4\lambda)}\cF[I_\lambda(f)](u),
\quad\text{$u\in\R^d$,}
$$
and let $h_\lambda=H_\lambda|_{Z}\in L_2(Z)$
 (thus $H_\lambda=E(h_\lambda)$).

Suppose that $k\in\N$.  Equation \eqref{E:1.3} implies that 
\begin{equation}\label{E:3fval1}
(2\pi)^d f(x_k)=\int_{ Z }\cF[f](u)e^{i\la x_k,u\ra}\,du=\big\la \cF[f],e^{-i\la x_k,\cdot\ra}\big\ra_{ Z }.
\end{equation}
On the other hand, equations \eqref{E:2.3.1} and \eqref{E:2.3.1a} assert that 
\begin{align}\label{E:3fval2}
(2\pi)^d f(x_k)
&=(2\pi)^d I_\lambda(f)(x_k) \\
&= \int_{\R^d}  \cF[I_\lambda(f)](u) e^{i\la x_k,u\ra} du \text{\ \ (by \eqref{E:1.3})}\notag\\
&= \int_{\R^d}    e^{-\|u\|_2^2/(4\lambda)}H_\lambda(u) e^{i\la x_k,u\ra}\, du\notag\\
&=\int_{ Z }    e^{-\|u\|_2^2/(4\lambda)}H_\lambda(u) e^{i\la x_k,u\ra} du\notag\\&\qquad +
\sum_{m=1}^\infty 
\int_{2^m Z \setminus 2^{m-1} Z }    e^{-\|u\|_2^2/(4\lambda)}H_\lambda(u) e^{i\la x_k,u\ra}\, du \notag\\
&=\int_{ Z }    e^{-\|u\|_2/^2(4\lambda)}H_\lambda(u) e^{i\la x_k,u\ra}\, du\notag\\ 
&\qquad+
\sum_{m=1}^\infty 2^{dm}
\int_{ Z \setminus 2^{-1} Z }    e^{-\|2^mv\|_2^2/(4\lambda)}H_\lambda(2^mv) e^{i\la x_k,2^mv \ra}\, dv \notag\\
&=\int_{ Z }    e^{-\|u\|_2^2/(4\lambda)}h_\lambda(u) e^{i\la x_k,u\ra}\, du\notag\\ 
&\qquad+
\sum_{m=1}^\infty 2^{dm}
\int_{ Z \setminus 2^{-1} Z }    e^{-\|2^mv\|_2^2/(4\lambda)} A_m(h_\lambda)(v)
 \overline{A_m (e^{-i\la x_k,\cdot \ra})} (v)\, dv\notag\\
&=\big\la   e^{-\|\cdot\|_2^2/(4\lambda)}h_\lambda ,e^{-i\la x_k,\cdot\ra} \big\ra_{ Z } \notag \\
&\qquad+\sum_{m=1}^\infty 2^{dm}\big\la  e^{-\|2^m(\cdot)\|_2^2/(4\lambda)} A_m(h_\lambda), A_m( e^{-i\la x_k,\cdot)\ra})\big\ra_{ Z }\notag\\ 
 &=\big\la   \cF[I_\lambda(f)] ,e^{-i\la x_k,\cdot\ra} \big\ra_{ Z } \notag\\
&\qquad
+\sum_{m=1}^\infty\big\la 2^{dm} A^*_m\big( e^{-\|2^m(\cdot)\|_2^2/4(\lambda)} A_m(h_\lambda)\big),  
e^{-i\la x_k,\cdot\ra}\big\ra_{ Z }\notag\\
&=\big\la   \cF[I_\lambda(f)]+\sum_{m=1}^\infty 2^{dm} A^*_m\big( e^{-\|2^m(\cdot)\|_2^2/4(\lambda)} A_m(h_\lambda)\big),  
e^{-i\la x_k,\cdot\ra}\big\ra_{ Z }.\notag
\end{align}
As $(e^{-i\la x_k,\cdot\ra}:k\in\N)$ is a Riesz basis for $L_2( Z )$ (in particular a complete system), 
equations
\eqref{E:3fval1} and \eqref{E:3fval2} lead to the identity
\begin{equation}\label{E:3.8}
\cF[f]=\cF[I_\lambda(f)]+\sum_{m=1}^\infty 2^{dm} A^*_m\big( e^{-\|2^m(\cdot)\|^2_2/4(\lambda)} A_m(h_\lambda)\big)\quad \text{ a.e. on $ Z $.}
\end{equation}
Suppose now that $h\in L_2( Z )$ and $m\in\N$.
We deduce from \eqref{E:3.7} that
\begin{align*}
\|2^{md}&A^*_m\big( e^{-\|2^m(\cdot)\|_2^2/(4\lambda)} A_m(h)\big) \|_{L_2( Z )}^2\\
&\le C^m R_b^4 2^{md}\|e^{-\|2^m(\cdot)\|_2^2/(4\lambda)} A_m(h)\big) \|_{L_2( Z )}^2\\
&\le C^m R_b^4 2^{md}\|e^{-2^{2m-2}\delta^2/(4\lambda)} A_m(h)\big) \|_{L_2( Z )}^2
\quad
\text{(since $\supp A_m(h)\!\subset\!  Z \setminus\frac12  Z\!\subset\!Z\setminus \frac\delta2 B_2 $)}\\
&\le \big(C^m R_b^4 \big)^2 e^{-2^{2m-2}\delta^2/(2\lambda)} \| h \|_{L_2( Z )}^2,
\end{align*}
whence
$$\|2^{md}A^*_m\big(e^{-\|2^m(\cdot)\|_2^2/(4\lambda)} A_m\big)\|_{L_2( Z )}
\le C^m R_b^4 e^{-2^{2m-2}\delta^2/(4\lambda)}.$$
Therefore the linear operator
$$\tau_\lambda: L_2( Z )\to L_2( Z ),\quad h\mapsto
\sum_{m\in\N}2^{md}A^*_m\big( e^{-\|2^m(\cdot)\|_2^2/(4\lambda)} A_m(h)\big)$$
is bounded. In fact, 
as
 there  are numbers $\lambda_0>0$ and $D$, which depend only on $C$ (which only depends on $d$),  such that
\begin{equation}\label{E:3.6c}
\sum_{m\in\N} C^{m} e^{-2^{2m-2}\delta^2/(4\lambda)}\le D e^{-\delta^2/(4\lambda)},\quad \lambda\in(0,\lambda_0],
\end{equation}
the operator norm of $\tau_\lambda$ obeys the following estimate:
 \begin{equation}\label{E:3.9}
 \|\tau_\lambda\|\le R_b^4 D e^{-\delta^2/(4\lambda)}  \text{ whenever $\lambda<\lambda_{0}$}.
 \end{equation}
 As the operator $\tau_\lambda$ is positive,  \eqref{E:3.8} yields
  $$\|\cF[f]\|_2\, \|h_\lambda\|_2\ge  
  \la \cF[f],h_\lambda\ra_{_Z}\ge\la e^{-\|\cdot\|_2^2/(4\lambda)} h_\lambda,h_\lambda\ra_{_Z}
  \ge e^{-1/(4\lambda)} \|h_\lambda\|_2^2.$$
Consequently,
\begin{equation}\label{E:3.9a}
\|h_\lambda\|_2\le e^{1/(4\lambda)} \|\cF[f]\|_2.\end{equation}
Thus, from \eqref{E:3.8}  and\eqref{E:3.9} we get
\begin{equation}\label{E:3.10}
\|\cF[I_\lambda(f)]|_{_{ Z }}\|_2
\le \|\cF[f]\|_2+ \|\tau_\lambda(h_\lambda)\|_2\le  \big(1+ R_b^4D e^{(1-\delta^2)/(4\lambda)}\big)\|\cF[f]\|_2.
\end{equation}
Our next step is to estimate $\|\cF[I_\lambda(f)]\,|_{_{\R^d\setminus  Z }}\|_2$.
Equation \eqref{E:3.1} implies that
\begin{align}\label{E:3.11a}
\|\cF[&I_\lambda(f)]|_{_{\R^d\setminus  Z }}\|_2^2\\
&=
\int_{\R^d\setminus  Z }e^{-\|u\|^2_2/(2\lambda)} |H_\lambda(u)|^2 \, du\notag\\
&=\sum_{m=1}^\infty \int_{2^m Z \setminus 2^{m-1} Z }e^{-\|u\|^2_2/(2\lambda)} |H_\lambda(u)|^2 \, du \notag\\
&=\sum_{m=1}^\infty 2^{dm}\int_{ Z \setminus 2^{-1} Z }e^{-2^{2m}\|v\|^2_2/(2\lambda)}|A_m(h_\lambda)(v)|^2 \, dv \notag\\
&\le \sum_{m=1}^\infty 2^{dm} e^{-2^{2m}\delta^2/(8\lambda)}\|A_m(h_\lambda)\|_2^2
\quad\text{(as $\supp A_m(h)\subset  Z \setminus2^{-1} Z $)}\notag\\
&\le R_b^4 \|h_\lambda\|^2_2\sum_{m=1}^\infty C^me^{-2^{2m}\delta^2/(8\lambda)} 
\quad\text{  \  (by \eqref{E:3.7})} 
\notag\\
&\le e^{1/(2\lambda)}  R_b^4  \|\cF[f]\|_2^2 \sum_{m=1}^\infty   e^{-2^{2m}\delta^2/(8\lambda)}C^m    
\quad \text{  \  (by \eqref{E:3.9a}).} \notag
\end{align}
By changing $\lambda_0$ and $D$, if need be, one obtains, as in \eqref{E:3.6c},
\begin{equation}\label{E:3.6c2}
\sum_{m=1}^\infty   e^{-2^{2m}\delta^2/(8\lambda)}C^m   \le De^{-\delta^2/(2\lambda)}, \quad \lambda\in(0,\lambda_0].
\end{equation}
Combining  \eqref{E:3.10} and \eqref{E:3.11a} proves our claim.

\medskip

\noindent
{\bf Step 2.}  Let $f\in PW_{\beta B_2}$.  There is a positive constant $D_2$ such that
\begin{equation}\label{E:3.13}
\|f-I_\lambda(f)\|_2 \leq D_2 e^{(\beta^2-3\delta^2+2)/(4\lambda)} \| f \|_2, 
\end{equation}
for all $0<\lambda<\lambda_0$.
\begin{rem}
  Note that \eqref{E:3.13} implies that 
$\lim_{\lambda\to0^+} I_\lambda(f)= f \text{\  in $L_2(\R^d)$.}$
\end{rem}
To prove \eqref{E:3.13} we define
 $\tilde\tau_\lambda=e^{1/(4\lambda)}\tau_\lambda$,  
\begin{align*}
&M_\lambda:L_2( Z )\to L_2( Z ), \quad h\mapsto e^{-(1-\|\cdot\|_2^2)/(4\lambda) }h, \text{ and }\\
&L_\lambda: L_2( Z ) \to  L_2( Z ), \quad h\mapsto R\circ\cF\circ I_\lambda\circ\cF^{-1}(h),
\end{align*}
where $R:L_2(\R^d)\to L_2( Z )$ is the restriction map.

\begin{prop}\label{P:3.2} 
The map $\Id+\tilde\tau_\lambda\circ M_\lambda$
is an invertible operator on $L_2( Z )$, and 
$(\Id+\tilde\tau_\lambda\circ M_\lambda)^{-1}=L_\lambda$.
\end{prop}

\begin{proof}
Let $h\in PW_{ Z }$.  From \eqref{E:3.8} we obtain (a.e. on $Z$)
\begin{align*}
\cF[h]&=\cF[I_\lambda(h)]+\tau_\lambda\big(e^{\|\cdot\|^2/(4\lambda)}\cF[I_\lambda(h)]|\big)\\
&=\cF[I_\lambda(h)]+\tilde\tau_\lambda\circ M_\lambda \big(\cF[I_\lambda(h)]\big)
 =\big(\Id+\tilde\tau_\lambda\circ M_\lambda \big) L_\lambda(\cF[h]).
\end{align*}
This implies that $\Id+\tilde\tau_\lambda\circ M_\lambda$ is surjective and is a left inverse of the bounded operator 
$L_\lambda$. Next we show that $\Id+\tilde\tau_\lambda\circ M_\lambda$ is also injective. To that end, 
let $(\Id+\tilde\tau_\lambda\circ M_\lambda)(h)=0 $ for some $h\in L_2( Z )$. Then
 \begin{align*}
  0&=\big\la (\Id\!+\!\tilde\tau_\lambda\circ M_\lambda)(h), M_\lambda(h)\big\ra_{ Z }
=\big\la h,M_\lambda(h)\big\ra_{ Z }\!+\!\big\la \tilde\tau_\lambda( M_\lambda(h)), M_\lambda(h)\big\ra_{ Z }
   \ge \big\la h,M_\lambda(h)\big\ra_{ Z }\ge 0,
 \end{align*}
the first inequality above being a consequence of the positivity of $\tilde\tau_\lambda$.  
Hence $\la h,M_\lambda(h)\ra_{ Z }=0$, which implies that $h=0$, because
$M_\lambda$ is a strictly positive operator.
The injectivity of  $\Id+\tilde\tau_\lambda\circ M_\lambda$ follows.
Thus  $\Id+\tilde\tau_\lambda\circ M_\lambda$ is invertible, and its
 inverse is $L_\lambda$. 
\end{proof}
Proposition \ref{P:3.2} provides the following identity on $ Z $:
$$
\cF[f]-\cF[I_\lambda (f)]=\big[\Id-(\Id+\tilde\tau_\lambda \circ M_\lambda)^{-1} \big](\cF[f])=
(\Id+\tilde\tau_\lambda\circ  M_\lambda)^{-1}\circ\tilde\tau_\lambda\circ M_\lambda(\cF[f]).
$$
If $f\in PW_{\beta B_2}$, then \eqref{E:3.9} and Step 1 provide
\begin{align}\label{E:3.16}
\| \cF[f]\!-\!\cF[I_\lambda (f)]|_{_{ Z }}\|_2&\le \|(\Id\!+\!\tilde\tau_\lambda\!\circ\!  M_\lambda)^{-1}\|\,\|\tilde\tau_\lambda\| \,\|M_\lambda(\cF[f])\|_2\\
                                 &\le D_1 e^{(1-\delta^2)/(4\lambda)} e^{1/(4\lambda) }R^2_b D e^{(-\delta^2)/(4\lambda)} \|M_\lambda(\cF[f])\|_2\notag\\
                  &=: D' e^{(1-2\delta^2)/(4\lambda)}   \big\|e^{\|\cdot\|^2/(4\lambda)}(\cF[f])\big\|_2\notag\\
				& \le D' e^{(\beta^2+1-2\delta^2)/(4\lambda)} \|\cF[f] \|.\notag
\end{align}

Now the first inequality in \eqref{E:3.11a} yields
\begin{align}\label{E:3.17}
 \big\|\cF[&I_\lambda (f)]|_{_{\R^d\setminus  Z }}\big\|_2^2\\
&\le \sum_{m=1}^\infty 2^{dm} e^{-2^{2m}\delta^2/(8\lambda)} \|A_m(h_\lambda)\|_2^2\notag\\
&\le  R_b^4\sum_{m=1}^\infty  C^m e^{-2^{2m}\delta^2/(8\lambda)} \big\|e^{\|\cdot\|^2_2/(4\lambda)} \cF[I_\lambda(f)]|{_{Z}}\big\|_2^2\quad
 \text{\ \ by \eqref{E:3.7}} \notag\\
 &\le  R_b^4 D e^{-\delta^2/(2\lambda)} \big\|e^{\|\cdot\|^2_2/(4\lambda)} \cF[I_\lambda(f)]|{_{Z}}\big\|_2^2\quad
 \text{\ \ by \eqref{E:3.6c2}} \notag\\
 &\le R_b^4 D\big[  \big\|e^{(\|\cdot\|^2_2-\delta^2)/(4\lambda)} (\cF[I_\lambda(f)]|{_{Z}}-\cF (f)\big)\big\|_2  + 
  \big\|e^{(\|\cdot\|^2_2-\delta^2)/(4\lambda)} \cF (f))\big\|_2 \big]^2 \notag\\
  &\le R_b^4 D\big[ e^{(1-\delta^2)/(4\lambda)} \big\| (\cF[I_\lambda(f)]|{_{Z}}-\cF (f)\big)\big\|_2  + 
  \big\|e^{(\|\cdot\|^2_2-\delta^2)/(4\lambda)} \cF (f))\big\|_2 \big]^2.\notag
 \end{align}
If we restrict to $f \in PW_{\beta B_2}$, then, by \eqref{E:3.16},
\begin{align*}
\big\|\cF[&I_\lambda (f)]|_{_{\R^d\setminus  Z }}\big\|_2^2 
 \le R_b^4 D\big[D' e^{(\beta^2+1-2\delta^2)/(4\lambda)}e^{(1-\delta^2)/(4\lambda)} + e^{(\beta^2-\delta^2)/(4\lambda)}\big]^2\|\cF (f))\|_2^2.
\end{align*}

Combining this  with \eqref{E:1.2}, and using \eqref{E:3.16} once again, we obtain 
the following estimate for some constant $D_2$:
$$ \Arrowvert f-I_\lambda f \Arrowvert_2 \leq D_2  e^{(\beta^2 -3\delta^2+2)/(4\lambda)} \Arrowvert f \Arrowvert_2. $$

\noindent{\bf Step 3.} Suppose that $f \in PW_{\beta B_2}$. There exist constants $\lambda_1\in(0,\lambda_0]$ and $D_3$  so that
$$\big|I_\lambda(f)(x)- f(x)\big|\le D_3 e^{(\beta^2-3\delta^2+2)/(4\lambda)}\|f\|_2,$$
for all $0<\lambda\le \lambda_1$ and  $x\in \R^d$. In particular
$\lim_{\lambda\to0^+} I_\lambda(f)= f$ uniformly on $\R^d$.

We first observe that, we can find, as before, numbers $\lambda_1\in(0,\lambda_0]$  and $D''>0$, such that
\begin{equation}\label{E:3.18}
 m^{1/2}( Z )R_b^2 \sum_{m=1}^\infty C^{m/2}  2^{dm/2}  e^{(4-2^{2m})/(16\lambda)}\le D'' \quad\text{whenever}\quad 0<\lambda\le \lambda_1. 
 \end{equation} 

Let $x\in\R^d$ and $f\in PW_{\beta B_2}$. We use \eqref{E:1.3} to write 
\begin{align}\label{E:3.18}
\big|I_\lambda&(f)(x)- f(x)\big|\\
&= \frac1{(2\pi)^d} \Bigg| \int_{ Z } \big[ \cF[I_\lambda(f)](u)   -  \cF[f](u) \big] e^{ixu} du
 + \int_{\R^d \setminus  Z }\cF[I_\lambda(f)](u)e^{ixu}  du\Bigg|\notag\\
&\le \frac1{(2\pi)^d}\big[ \|\cF[I_\lambda(f)]|_{_{ Z }}- \cF[f]\|_1+\big\|\cF[I_\lambda(f)]|_{_{\R^d\setminus  Z }}\|_1\big].\notag
\end{align}
From the BCS  inequality and \eqref{E:3.16}  we deduce that  $\lim_{\lambda\to 0^+}\|\cF[I_\lambda(f)]|_{_{ Z }}- \cF[f]\|_1=0$,
and an argument similar to that in \eqref{E:3.11a} yields
\begin{align*}
 \big\|\cF&[I_\lambda(f)]|_{_{\R^d\setminus  Z }}\big\|_1\\
&=\sum_{m=1}^\infty \int_{2^m  Z \setminus 2^{m-1}  Z  } e^{-\|u\|^2_2/(4\lambda)} |H_\lambda(u)|du\\
&=\sum_{m=1}^\infty 2^{dm}\int_{  Z \setminus 2^{-1}  Z  } e^{-2^{2m}\|v\|^2_2/(4\lambda)}|A_m(h_\lambda)(v)| dv\\
&\le m^{1/2}( Z )\sum_{m=1}^\infty 2^{dm} \|e^{-2^{2m}\|\cdot\|^2_2/(4\lambda)}A_m(h_\lambda)\|_2
\qquad\text{\ (by the BCS inequality)}\\
&\le  m^{1/2}( Z )R_b^2 \sum_{m=1}^\infty  C^{m/2} 2^{dm/2}  e^{-2^{2m}\delta^2/(16\lambda)} \|h_\lambda\|_2\\
&\qquad\qquad\Big(\text{\ by \eqref{E:3.7} and since $\supp(A_m(h))\!\subset\!Z\setminus { \frac12}Z$}\Big)\\
&= D'' \big\|e^{(\|\cdot\|_2^2-\delta^2)/(4\lambda)}\cF[I_\lambda(f)]|_{ Z }\big\|_2   \quad\text{(by \eqref{E:3.18})}\\
&\leq D''  e^{(1-\delta^2)/(4\lambda)}\big[\|\cF[I_\lambda(f)]|_Z-\cF[f]\|_2 +\|\cF[f]\|_2\big]\\
&\leq D'D''e^{(\beta^2+2-3\delta^2)/(4\lambda)}\|\cF[f]\|. \quad\text{(by \eqref{E:3.16})}
\end{align*} 
This concludes the proof.
\end{proof}

\centerline{\bf Acknowledgments\/}

\noindent We thank Yuri~Lyubarskii, Joaquim~Ortega-Cerd\`a, Grigoris~Paouris, and Kristian~Seip for generously sharing with us their time and expertise.

\end{document}